\documentclass[12pt]{article}
\input{amssym}
\input{epsf}

\begin{document}
\title{{Preliminary group classification for generalized inviscid Burger's equation}}
\date{}

\author{
%M.A. Jafarizadeh \\ {\small{{\it Physics department, Tabriz University, Tabriz,
%Iran.}}}\\ {\small jafarizadeh@tabrizu.ac.ir}\\[2mm] \and
A. Mahdipour--Shirayeh\thanks{Corresponding author.} \\
{\small{{\it Department of Mathematics, Semnan University, Semnan,
I.R.Iran}}}\\ {\small ali.mahdipour@gmail.com}\\
%\and T. \"Ozer \\
%{\small{{\it Istanbul Technical University, Faculty of Civil
%Engineering, Division of Mechanics,}}}\\ {\small{\it No. 36234469
%Maslak, Istanbul, Turkey}}\\{\small tozer@itu.edu.tr}
%\and F. Rezvan \\
%{\small{{\it Faculty of Mathematical Sciences, Tabriz University,
%Tabriz, I.R.Iran}}}\\ {\small frezvan@tabrizu.ac.ir}
}
\maketitle
\begin{abstract}
Preliminary group classification for a class of generalized
inviscid Burger's equations in the general form $u_t + g(x,u)u_x =
f(x,u)$ is given and additional equivalence transformations are
found. Adduced results complete and essentially generalize recent
works on the subject
%[M. Nadjafikhah, Lie Symmetries of Inviscid
%Burgers' Equation, Adv. appl. Clifford alg., DOI
%10.1007/s00006-003-0000; N. Smaoui and M. Mekkaoui, The
%generalized Burgers equation with and without a time delay, J.
%Applied Math. Stochastic Anal. Vol. 1 (2004), pp. 73-96.]
. A number of new interesting nonlinear invariant models which
have non--trivial invariance algebras are obtained. The result of
the work is a wide class of equations summarized in table form.
\end{abstract}
%
%\medskip \noindent {\bf PACS numbers:} 47.35.+i, 47.40.-x
\noindent {\bf Key words}: Burgers' equation, Equivalence
transformations, Lie symmetries.
%
%%%%%%%%%%%%%%%%%%%%%%%%%%%%%%%%%%%%%%%%%%%%%%%%%%%%%%%%%%%
\section{Introduction}
The theory of Lie groups has greatly influenced diverse branches
of mathematics and physics. The main tool of the theory, Sophus
Lie's infinitesimal method \cite{Lie}, establishes connection
between continuous transformation groups and algebras of their
infinitesimal generators. The method leads to many techniques of
great significance in studying the group-invariant solutions and
conservation laws of differential equations \cite{Ib2,Ol,Ov}. The
approach to the classification of partial differential equations
which we propound is, in fact, a synthesis of Lie's infinitesimal
method, the use of equivalence transformations and the theory of
classification of abstract finite-dimensional Lie algebras.

Historically, Burgers' equation is a fundamental partial
differential equation from fluid mechanics \cite{Ib4}. There have
been a number of papers that have contributed to the studies of
Lie groups of transformations of various classes of Burgers'
equation, such as modeling of gas dynamics and traffic flow
\cite{PZ}.

The classical form of Burgers' equation is
\begin{eqnarray}
u _t + u \,u_x = \nu\, u_{xx}, \label{eq:1}
\end{eqnarray}
for viscosity coefficient $\nu$, which describes the evolution of
the field $u=u(t,x)$ under nonlinear advection and linear
dissipation [3, 10]. When $\nu = 0$, Burgers' equation becomes the
inviscid Burgers' equation:
\begin{eqnarray}
u _t + u \,u_x = 0. \label{eq:1-1}
\end{eqnarray}
which is a prototype for equations for which the solution can
develop discontinuities (shock waves). The generalized
non--homogeneous inviscid Burgers' equation had introduced in
\cite{SM} has the form
\begin{eqnarray}
u_t + g(u)\,u_x = f(u). \label{eq:1-2}
\end{eqnarray}

In this study, we extend equation (\ref{eq:1-2}) and consider the
following generalized non--homogeneous inviscid Burgers' equation
(NIB) in the general form
\begin{eqnarray}
\hspace{-0.8cm}{\rm NIB}:\hspace{0.5cm}  u _t + g(x,u) \,u_x =
f(x,u), \label{eq:2}
\end{eqnarray}
where $f\neq 0$ and $g\neq 0$ are sufficiently smooth nonconstant
functions of variables $x,u$.

The main results of the present paper are point symmetry and
equivalence classification of NIB equations leading to a number of
new interesting nonlinear invariant models associated to
non--trivial invariance algebras. A complete list of these models
are given for a finite--dimensional subalgebra of
infinite--dimensional equivalence algebra derived for NIB
equations. Reaching to these goals we follow an algorithm
explained and performed in references \cite{Ib3,Ov}. Our method is
similar to the way of \cite{Ma} for the nonlinear heat
conductivity equation $u_t=\left[E(x,u)u_x\right]_x + H(x,u)$.

Throughout this paper we assume that each index of a function
implies the derivation of the function with respect to it, unless
specially stated otherwise.
%
%%%%%%%%%%%%%%%%%%%%%%%%%%%%%%%%%%%%%%%%%%%%%%%%%%%%%%%%%%%%%%
%
\section{Principal Lie algebra}

The symmetry approach to the classification of admissible partial
differential equations relies heavily upon the availability of a
constructive way of describing transformation groups leaving
invariant the form of a given partial differential equation. This
is done via the well-known infinitesimal method developed by
Sophus Lie \cite{Ol,Ol2,Ov}. Given a partial differential
equation, the problem of investigating its maximal (in some sense)
Lie invariance group reduces to solving an over--determined system
of linear partial differential equations, called the {\it
determining equations}.

Considering the total space $E=X\times U$ with local coordinate
$(t,x,u)$ of independent variables $(t,x)\in X$ and dependent
variable $u\in U$, the solution space of Eq.~(\ref{eq:2}), is a
subvariety $S_{\Delta}\subset J^1({\Bbb R}^2,{\Bbb R})$ of the
first order jet bundle of 2-dimensional submanifolds of $E$. The
one--parameter Lie group of infinitesimal transformations on $E$
is as follows
\begin{eqnarray}
&& \tilde{t}=t + s\,\xi(t,x,u) + O(s^2),   \nonumber\\
&& \tilde{x}=x + s\,\tau(t,x,u) + O(s^2),   \label{eq:3}\\
&& \tilde{u}=u + s\,\varphi(t,x,u) + O(s^2),   \nonumber ,
\end{eqnarray}
where $s$ is the group parameter and $\xi, \tau$ and $\varphi$ are
the (sufficiently smooth) infinitesimals of the transformations
for the independent and dependent variables, resp. Thus the
corresponding infinitesimal generators have the following general
form
\begin{eqnarray}
{\bf v}=\xi(t,x,u)\,\frac{\partial}{\partial t}+
\tau(t,x,u)\,\frac{\partial}{\partial
x}+\varphi(t,x,u)\,\frac{\partial}{\partial u}.\label{eq:4}
\end{eqnarray}
According to the invariance condition (see e.g. \cite{Ib}, Theorem
2.36), ${\bf v}$ is a point transformation of Eq.~(\ref{eq:2}) if
and only if
\begin{eqnarray}
{\bf v}^{(1)}\Big[u _t + g\,u_x - f
\Big]\Big|_{Eq.~(\ref{eq:2})}=0. \label{eq:5}
\end{eqnarray}
In the latter relation, ${\bf v}^{(1)}$ is the first order
prolongation of vector field ${\bf v}$:
\begin{eqnarray}
{\bf v}^{(1)}={\bf v} +
\varphi^t(t,x,u,u_t,u_x,u_{tt},u_{tx},u_{xx})\,\frac{\partial}{\partial
u_t} +
\varphi^x(t,x,u,u_t,u_x,u_{tt},u_{tx},u_{xx})\,\frac{\partial}{\partial
u_x},\label{eq:4}
\end{eqnarray}
in which
\begin{eqnarray}
\varphi^t={\Bbb D}_t({\cal Q}) + \xi\,u_{tt}+\tau\,u_{tx},\\
\varphi^x={\Bbb D}_x({\cal Q}) + \xi\,u_{tx}+\tau\,u_{xx},
\end{eqnarray}
and ${\cal Q}=\varphi - \xi\,u_t-\eta\,u_x$ is the characteristic
of vector field ${\bf v}$ and the operators $${\Bbb
D}_J=\partial_J + u_{t,J}\partial u_t + u_{x,J}\partial u_x +
\cdots$$ (for $J=t,x$) are total derivatives with respect to $t$
and $x$. Applying ${\bf v}^{(1)}$ on Eq.~(\ref{eq:2}) and
introducing the relation $u_t= f - g\,u_x$ to eliminate $u_t$, we
obtain the following determining equation:
\begin{eqnarray}
&& \hspace{-1.5cm} -\tau\,g\,g_x\,u_x + \tau\,f\,g_x - \tau\,f_x +
\varphi\,g_u\,u_x - \varphi\,f_u + \varphi_t + f\,\varphi_u
-g\,\varphi_u\,u_x + g\,\xi_t\,u_x - f\,\xi_t   \nonumber\\[-3mm]
&& \label{eq:6}\\[-2.5mm]
&& \hspace{-1.5cm} - \tau_t\,u_x - f^2\,\xi_u + f\,g\,\xi_u\,u_x -
\tau_u\,f\,u_x + g\,\varphi_x + g\,\varphi_u\,u_x -g\,\tau_x\,u_x
+ g^2\,\xi_x\,u_x - f\,g\,\xi_x =0. \nonumber
\end{eqnarray}
In the last relation, since $\xi, \tau$ and $\varphi$ do not
depend to variable $u_x$ and its powers, so Eq.~(\ref{eq:6}) is
hold if and only if the following equations are fulfilled
\begin{eqnarray}
&& -\tau\,g\,g_x + \varphi\,g_u + \xi_t\,g - \tau_t + f\,g\,\xi_u - \tau_u\,f - g\,\tau_x+g^2\xi_x=0,    \nonumber\\[-3mm]
&& \label{eq:7}\\[-2.5mm]
&& \varphi_t -\tau\,f_x - \varphi\,f_u + \varphi_u\,f - \xi_t\,f -
\xi_u\,f^2 + g\,\varphi_x - f\,g\,\xi_x=0. \nonumber
\end{eqnarray}
But using the fact that $\{1,g,g^2\}$ and $\{1,f,f^2\}$ for
nonconstant functions $f,g$ are independent sets (see Lemma 1 of
\cite{Na}), one can divide Eqs.~(\ref{eq:7}) to the below
over--determined system
\begin{eqnarray}
&& \tau_u=0,\hspace{0.5cm}  \xi_u=0,\hspace{0.6cm}
\xi_t -\tau_x-\tau\,g_x =0,\hspace{0.9cm}  \varphi\,g_u - \tau_t=0,    \nonumber\\[-3mm]
&& \label{eq:8}\\[-2.5mm]
&& \xi_x=0,\hspace{0.5cm} \varphi_x=0,\hspace{0.5cm}
\phi_t-\tau\,f_x-\varphi\,f_u =0, \hspace{0.5cm} \varphi_u -
\xi_t= 0. \nonumber
\end{eqnarray}
General solution to this system results in the final form of point
infinitesimal generators:

\paragraph{Theorem 1.} {\em Infinitesimal generator of every one--parameter Lie group of point
symmetries of NIB equation is
\begin{eqnarray}
{\bf v}=\frac{\partial}{\partial t}. \label{eq:8-1}
\end{eqnarray}
Furthermore, every NIB equation in the form Eq.~(\ref{eq:2}) can
not be reduced into an inhomogeneous form of a linear equation.}\\

\noindent{\it Proof:} The second assertion is a simple result of
p. 209 of \cite{Ol}. \hfill\ $\diamondsuit$\\

Let we change the conditions on $f,g$ and permit some of them be
zero: Assume that $f=0$ and $g=g(u)$. In this case, NIB equation
reduced to the the equation $IBE:~ u_t+g(u)\,u_x=0$ which has
recently studied in \cite{Na}. Applying the conditions on $f, g$
on determining equation (\ref{eq:6}) and the same computation as
above show that the general form of infinitesimal, represented in
Theorem 1 of \cite{Na} is not correct. The correct form of
infinitesimal generators of IBE is as follows
\begin{eqnarray}
{\bf v} = [F^2(u)\,t+F^1(u)]\frac{\partial}{\partial t}
+[F^4(u)\,g_u\,t + F^2(u)\,x + F^3(u)]\frac{\partial}{\partial x}
+ F^4(u)\,\frac{\partial}{\partial u},
\end{eqnarray}
where $F^i$~s are arbitrary smooth functions of $u$. One may
divides ${\bf v}$ to the following vector fields
\begin{eqnarray}
&& {\bf v}_1:={\bf v}^1_{F^1} = F^1(u)\,\frac{\partial}{\partial
t},\hspace{1cm}
{\bf v}_2:={\bf v}^2_{F^2}= F^2(u)\Big[t\,\frac{\partial}{\partial
t}
+x\frac{\partial}{\partial x}\Big],   \nonumber\\[-3mm]
&& \label{eq:9}\\[-2.5mm]
&& {\bf v}_3:={\bf v}^3_{F^3} = F^3(u)\,\frac{\partial}{\partial
x},\hspace{1cm}
{\bf v}_4:={\bf v}^4_{F^4} =
F^4(u)\Big[g_u\,t\,\frac{\partial}{\partial x} +
\frac{\partial}{\partial u}\Big]. \nonumber
\end{eqnarray}
The commutative table of these vector fields is given in Table 1.
\begin{table}
\centering{\caption{{\small The commutators table for point
symmetry algebra of IBE.}}}\label{table:1}\vspace{-0.5cm}
\begin{eqnarray*}\begin{array}{l|l l l l}
\hline\hline
      &\hspace{0.1cm} {\bf v}_1   &\hspace{1cm} {\bf v}_2              &\hspace{1cm} {\bf v}_3      &\hspace{1cm} {\bf v}_4  \\ \hline
  {\bf v}_1 &\hspace{0.1cm} 0     &\hspace{1cm} {\bf v}^1_{F^1F^2}     &\hspace{1cm} 0    &\hspace{1cm} -{\bf v}^1_{F^1_uF^4}+{\bf v}^3_{F^1_uF^4 g_u} \\
  {\bf v}_2 &\hspace{0.1cm} -{\bf v}^1_{F^1F^2}   &\hspace{1cm} 0      &\hspace{1cm} -{\bf v}^3_{F^2F^3}   &\hspace{1cm} -{\bf v}^2_{F^2_uF^4} \\
  {\bf v}_3 &\hspace{0.1cm} 0        &\hspace{1cm} {\bf v}^3_{F^2F^3}  &\hspace{1cm} 0      &\hspace{1cm} {\bf v}^3_{F^3_uF^4} \\
  {\bf v}_4 &\hspace{0.1cm} {\bf v}^1_{F^1_uF^4}-{\bf v}^3_{F^1_uF^4 g_u}    &\hspace{1cm}{\bf v}^2_{F^2_uF^4}  &\hspace{1cm} -{\bf v}^3_{F^3_uF^4} &\hspace{1cm} 0 \\
  \hline\hline
 \end{array}\end{eqnarray*}
\end{table}

Also every projective infinitesimal generator introduced in
Theorem 2 of \cite{Na} has the following correct form
\begin{eqnarray}
{\bf v} = [c_1\,t+F^1(u)]\frac{\partial}{\partial t} +[c_3\,t +
c_1\,x + c_2]\frac{\partial}{\partial x} +
\frac{c_3}{g_u}\,\frac{\partial}{\partial u},
\end{eqnarray}
for $g_u\neq0$. For $g(u)\equiv \verb"const."$, it has the
following modified form
\begin{eqnarray}
{\bf v} = [c_1\,t+F^1(u)]\frac{\partial}{\partial t} +[ c_1\,x +
c_2]\frac{\partial}{\partial x} + F^2(u)\,\frac{\partial}{\partial
u}.
\end{eqnarray}
A generalized case of IBE occurs when $f=0$ and $g=g(x,u)$. Thus
we deal with the homogeneous inviscid Burgers' equation. The
similar computations show that in this case, the general form of
point infinitesimal generators is
\begin{eqnarray}
{\bf v} = [F^2(u)\,t+F^1(u)]\frac{\partial}{\partial t}
+e^{-g}\Big[F^2\int\,e^g\,dx + F^3(u)\Big]\frac{\partial}{\partial
x}.
\end{eqnarray}
%
%%%%%%%%%%%%%%%%%%%%%%%%%%%%%%%%%%%%%%%%%%%%%%%%%%%%%%%%%%%%%%
%
\section{Equivalence transformations}
It is well--known that an equivalence transformation is a
non--degenerate change of the variables $t,x,u$ taking any
equation of the form (\ref{eq:2}) into an equation of the same
form, generally speaking, with different $f(x,u)$ and $g(x,u)$:
The set of all equivalence transformations forms an equivalence
group ${\cal E}$: We shall find a continuous subgroup ${\cal E}_c$
of it making use of the infinitesimal method \cite{Ib3,Ov}.

To find the group ${\cal E}_c$ we need to determine those
infinitesimal generators
\begin{eqnarray}
Y = \xi(t,x,u)\frac{\partial}{\partial
t}+\tau(t,x,u)\frac{\partial}{\partial
x}+\varphi(t,x,u)\frac{\partial}{\partial
u}+\chi(t,x,u,f,g)\frac{\partial}{\partial
f}+\eta(t,x,u,f,g)\frac{\partial}{\partial g},\label{eq:9-1}
\end{eqnarray}
from the invariance conditions of Eq.~(\ref{eq:2}) as the
following system
\begin{eqnarray}\hspace{-1.3cm}\left\{\begin{array}{l}
u_t + g(x,u)\,u_x=f(x,u)\\
f_t=g_t=0.
\end{array}\right.\label{eq:10}
\end{eqnarray}
Here, $u$ and $f , g$ are considered as differential variables:
$u$ on the space $(t, x)$ and $f$ , $g$ on the extended space $(t,
x, u)$. The coordinates $\xi,\tau, \varphi$ are sought as
functions of $t, x, u$ while the coordinates $\chi, \eta$ are
sought as functions of $t, x, u, f , g$. The invariance conditions
of the system (\ref{eq:10}) are
\begin{eqnarray}\left\{\begin{array}{l}
\tilde{Y}[u_t + g(x,u)\,u_x-f(x,u)]=0,\\
\tilde{Y}[f_t]=\tilde{Y}[g_t]=0,
\end{array}\right.\label{eq:11}
\end{eqnarray}
where $\tilde{Y}$ is the prolongation of (\ref{eq:9-1}) to the
first order jet space of differential variables $t,x,u,f,g$ which
one may represent as follows
\begin{eqnarray}
\tilde{Y} = Y + \varphi^t\,\frac{\partial}{\partial u_t}+
\varphi^x \,\frac{\partial}{\partial u_x}
+\chi^t\,\frac{\partial}{\partial
f_t}+\eta^t\,\frac{\partial}{\partial g_t}.\label{eq:12}
\end{eqnarray}
In this relation $\varphi^t$ and $\varphi^x$ are the same with
those introduced in section 2 and
\begin{eqnarray}
&&\hspace{-2cm} \chi^t = \widetilde{{\mathcal D}}_t(\chi) -
f_t\,\widetilde{{\mathcal D}}_t(\xi)- f_x\,\widetilde{{\mathcal
D}}_t(\tau) - f_u\,\widetilde{{\mathcal D}}_t(\varphi) =
\widetilde{{\mathcal D}}_t(\chi) - f_x\,\widetilde{{\mathcal
D}}_t(\tau) - f_u\,\widetilde{{\mathcal D}}_t(\varphi),  \label{eq:13}\\
&&\hspace{-2cm}  \eta^t = \widetilde{{\mathcal D}}_t(\eta) -
g_t\,\widetilde{{\mathcal D}}_t(\xi)- g_x\,\widetilde{{\mathcal
D}}_u(\tau) - g_u\,\widetilde{{\mathcal D}}_u(\varphi)
=\widetilde{{\mathcal D}}_t(\eta) - g_x\,\widetilde{{\mathcal
D}}_t(\tau) - g_u\,\widetilde{{\mathcal D}}_u(\varphi),
\label{eq:14}
\end{eqnarray}
where we consider
\begin{eqnarray}
\widetilde{{\mathcal D}}_t := \displaystyle{\frac{\partial
}{\partial t}} + f_t\,\displaystyle{\frac{\partial }{\partial f}}
+ g_t\,\displaystyle{\frac{\partial}{\partial g}} =
\displaystyle{\frac{\partial }{\partial t}}. \label{eq:15-0}
\end{eqnarray}
Substituting relations (\ref{eq:13})-(\ref{eq:15-0}) in
(\ref{eq:12}) and then applying it on the last two equations of
(\ref{eq:11}) we find that
\begin{eqnarray}
&& \chi_t-f_x\,\tau_t-f_u\,\varphi_t=0,\\
&& \eta_t - g_x\,\tau_t - g_u\,\varphi_t=0. \label{eq:15-1}
\end{eqnarray}
Since the latter equations are hold for every $f$ and $g$, so we
lead to the following facts
\begin{eqnarray}
\tau_t=0, \hspace{1cm}   \varphi_t=0,  \hspace{1cm}  \chi_t=0,
\hspace{1cm} \eta_t=0. \label{eq:16}
\end{eqnarray}
By effecting (\ref{eq:12}) on the first equation of (\ref{eq:11}),
we have
\begin{eqnarray}
-\chi + \eta\,u_x+\varphi^t + g\,\varphi^x = 0. \label{eq:17}
\end{eqnarray}
Substituting $\varphi^t, \varphi^x$, introducing $u_t=-g\,u_x+f$
to eliminate $u_t$ and using the fact that in the derived relation
$u_x$ and its powers are free variables, finally we tend to the
below system
\begin{eqnarray}
&& f\,g\,\xi_u + g\,\xi_t - f\,\tau_u -g\,\tau_x + g^2\,\xi_x +
\eta=0,    \label{eq:18}\\
&& -\chi + f\,\varphi_u - f\,\xi_t -f^2\,\xi_u + g\,\varphi_x
-f\,g\,\xi_x = 0. \label{eq:19}
\end{eqnarray}
The general solution to Eqs.~(\ref{eq:12}), (\ref{eq:12}) and
(\ref{eq:12}) is
\begin{eqnarray}
&& \xi=\xi(t,x,u), \hspace{1cm}   \eta=\tau(x,u), \hspace{1cm}    \varphi=\varphi(x,u)    \label{eq:20}\\
&& \chi=g\,\varphi_x - f\,g\,\xi_x - f\,\xi_t + f\,\varphi_u -
f^2\,\xi_u, \label{eq:21}\\
&& \eta=f\,\tau_u + g\,\tau_x - g\,\xi_t - g^2\xi_x - f\,g\,\xi_u.
\label{eq:22}
\end{eqnarray}
After utilizing these relations in $Y$, one may divide $Y$ to the
following vector fields
\begin{eqnarray}
&&\hspace{-1cm} Y_1:=Y^1_{\xi}=\xi\,\frac{\partial}{\partial t}
-f[g\,\xi_x+\xi_t+f\,\xi_u]\,\frac{\partial}{\partial f}
-g[f\,\xi_u + \xi_t + g\,\xi_x]\,\frac{\partial}{\partial g},    \nonumber\\
&&\hspace{-1cm}  Y_2:=Y^2_{\tau}=\tau\,\frac{\partial}{\partial x} +
[g\,\tau_x + f\,\tau_u]\,\frac{\partial}{\partial g},  \label{eq:24}\\
&&\hspace{-1cm}
Y_3:=Y^3_{\varphi}=\varphi\,\frac{\partial}{\partial u} +
[g\,\varphi_x + f\,\varphi_u]\,\frac{\partial}{\partial f}.
\nonumber \label{eq:25}
\end{eqnarray}
The commutators table of $Y_i$~s is given in Table 2.
\begin{table}
\centering{\caption{{\small The commutators table for equivalence
symmetry algebra of NIB.}}}\label{table:2}\vspace{-0.5cm}
\begin{eqnarray*}\begin{array}{l|l l l}
\hline\hline
      &\hspace{0.1cm} Y_1   &\hspace{1cm} Y_2              &\hspace{1cm} Y_3       \\ \hline
  Y_1 &\hspace{0.1cm} 0     &\hspace{1cm} -Y^1_{\tau\,\xi_x}     &\hspace{1cm} -Y^1_{\varphi\,\xi_u}    \\
  Y_2 &\hspace{0.1cm} Y^1_{\tau\,\xi_x}   &\hspace{1cm} 0      &\hspace{1cm} Y^3_{\tau\,\varphi_x}-Y^2_{\varphi\,\tau_u}    \\
  Y_3 &\hspace{0.1cm} Y^1_{\varphi\,\xi_u}        &\hspace{1cm} Y^2_{\varphi\,\tau_u}-Y^3_{\tau\,\varphi_x}  &\hspace{1cm} 0      \\
  \hline\hline
 \end{array}\end{eqnarray*}
\end{table}
%
%
%%%%%%%%%%%%%%%%%%%%%%%%%%%%%%%%%%%%%%%%%%%%%%%%%%%%%%%%%%%%%%
%
\section{Preliminary group classification}
In diverse applications of symmetry analysis, the equivalence
symmetry is handled to extend the principal Lie algebra to the
equivalence algebra ${\cal L}_{\cal F}$ when some further
equations under considerations are taken. These extensions are
called ${\cal F}$--extensions of the principal Lie algebra. The
classification of all non--equivalent equations (with respect to a
given equivalence group $G_{\cal F}$) admitting ${\cal
F}$--extensions of the principal Lie algebra is called a {\it
preliminary group classification}.

In the general, $G_{\cal F}$ is not necessarily the largest
equivalence group. In the following, we consider a subgroup of the
group of all equivalence transformations that has a
finite-dimensional subalgebra (desirably as large as possible) of
an infinite-dimensional algebra with basis (\ref{eq:24}) and use
it for a preliminary group classification. We select the
subalgebra ${\cal L}_{10}$, generated by the following vector
fields
\begin{eqnarray}
&&\hspace{-1cm} X_1=\frac{\partial}{\partial t}, \hspace{1cm}
X_2=\frac{\partial}{\partial x}, \hspace{1cm}
X_3=\frac{\partial}{\partial u},   \nonumber\\
&&\hspace{-1cm} X_4= t\,\frac{\partial}{\partial
t}-f\,\frac{\partial}{\partial f}-g\,\frac{\partial}{\partial
g},\hspace{1cm}  X_5= x\,\frac{\partial}{\partial
t}-f\,g\,\frac{\partial}{\partial f}-g^2\,\frac{\partial}{\partial
g},  \nonumber\\
&&\hspace{-1cm} X_6= u\,\frac{\partial}{\partial
t}-f^2\,\frac{\partial}{\partial f}-f\,g\,\frac{\partial}{\partial
g},\hspace{1cm}  X_7= x\,\frac{\partial}{\partial
x}+g\,\frac{\partial}{\partial
g},    \label{eq:26}\\
&&\hspace{-1cm} X_8= u\,\frac{\partial}{\partial
x}+f\,\frac{\partial}{\partial g},\hspace{1cm} X_9=
x\,\frac{\partial}{\partial
u}+g\,\frac{\partial}{\partial f},  \nonumber\\
&&\hspace{-1cm} X_{10}= u\,\frac{\partial}{\partial
u}+f\,\frac{\partial}{\partial f},   \nonumber
\end{eqnarray}
The commutator and adjoint representations of ${\cal L}_{10}$ are
listed in Tables 3 and 4.

It is well known that the problem of the construction of the
optimal system of solutions is equivalent to that of the
construction of the optimal system of subalgebras
\cite{Ol,Ol2,Ov}. Here, we determine a list (an optimal system) of
conjugacy inequivalent subalgebras with the property that any
other subalgebra is equivalent to a unique member of the list
under some element of the adjoint representation i.e.
$\overline{{\goth h}}\,{\rm Ad}(g)\,{\goth h}$ for some $g$ of a
considered Lie group. The adjoint action is given by the Lie
series
\begin{eqnarray}
{\rm Ad}(\exp(s\,Y_i))Y_j
=Y_j-s\,[Y_i,Y_j]+\frac{s^2}{2}\,[Y_i,[Y_i,Y_j]]-\cdots,
\end{eqnarray}
Then we will deal with the construction of the optimal system of
subalgebras of ${\cal L}_{10}$.
\begin{table}
\centering{\caption{{\small Commutators table for ${\cal
L}_{10}$}}}\label{table:3}\vspace{-0.5cm}
\begin{eqnarray*}\begin{array}{l|c c c c c c c c c c}
\hline\hline
      &\hspace{0.1cm} X_1 &\hspace{0.25cm} X_2&\hspace{0.25cm}X_3 &\hspace{0.25cm}X_4 &\hspace{0.25cm}X_5 &\hspace{0.25cm} X_6&\hspace{0.25cm}X_7 &\hspace{0.25cm} X_8      &\hspace{0.25cm}X_9       &\hspace{0.25cm} X_{10}  \\ \hline
  X_1 &\hspace{0.1cm}   0 &\hspace{0.25cm} 0  &\hspace{0.25cm}0   &\hspace{0.25cm}X_1 &\hspace{0.25cm}0   &\hspace{0.25cm} 0  &\hspace{0.25cm}0   &\hspace{0.25cm} 0        &\hspace{0.25cm}0         &\hspace{0.25cm} 0  \\
  X_2 &\hspace{0.1cm}   0 &\hspace{0.25cm} 0  &\hspace{0.25cm}0   &\hspace{0.25cm} 0  &\hspace{0.25cm}X_1 &\hspace{0.25cm} 0  &\hspace{0.25cm}X_2 &\hspace{0.25cm} 0        &\hspace{0.25cm}X_3       &\hspace{0.25cm} 0  \\
  X_3 &\hspace{0.1cm}   0 &\hspace{0.25cm} 0  &\hspace{0.25cm}0   &\hspace{0.25cm} 0  &\hspace{0.25cm}0   &\hspace{0.25cm} X_1&\hspace{0.25cm}0   &\hspace{0.25cm} X_2      &\hspace{0.25cm}0         &\hspace{0.25cm} X_3  \\
  X_4 &\hspace{0.1cm}-X_1 &\hspace{0.25cm} 0  &\hspace{0.25cm}0   &\hspace{0.25cm} 0  &\hspace{0.25cm}-X_4&\hspace{0.25cm}-X_5&\hspace{0.25cm}0   &\hspace{0.25cm} 0        &\hspace{0.25cm}0         &\hspace{0.25cm} 0  \\
  X_5 &\hspace{0.1cm}   0 &\hspace{0.25cm}-X_1&\hspace{0.25cm}0   &\hspace{0.25cm}X_4 &\hspace{0.25cm}0   &\hspace{0.25cm} 0  &\hspace{0.25cm}-X_4&\hspace{0.25cm}-X_5      &\hspace{0.25cm}0         &\hspace{0.25cm} 0  \\
  X_6 &\hspace{0.1cm}   0 &\hspace{0.25cm} 0  &\hspace{0.25cm}-X_1&\hspace{0.25cm}X_5 &\hspace{0.25cm}0   &\hspace{0.25cm} 0  &\hspace{0.25cm}0   &\hspace{0.25cm} 0        &\hspace{0.25cm}-X_4      &\hspace{0.25cm}-X_5  \\
  X_7 &\hspace{0.1cm}   0 &\hspace{0.25cm}-X_2&\hspace{0.25cm}0   &\hspace{0.25cm} 0  &\hspace{0.25cm}X_4 &\hspace{0.25cm} 0  &\hspace{0.25cm}0   &\hspace{0.25cm}-X_7      &\hspace{0.25cm}X_9       &\hspace{0.25cm} 0  \\
  X_8 &\hspace{0.1cm}   0 &\hspace{0.25cm} 0  &\hspace{0.25cm}-X_2&\hspace{0.25cm} 0  &\hspace{0.25cm}X_5 &\hspace{0.25cm} 0  &\hspace{0.25cm}X_7 &\hspace{0.25cm} 0        &\hspace{0.25cm}X_{10}-X_6&\hspace{0.25cm} X_7  \\
  X_9 &\hspace{0.1cm}   0 &\hspace{0.25cm}-X_3&\hspace{0.25cm}0   &\hspace{0.25cm} 0  &\hspace{0.25cm}0   &\hspace{0.25cm} X_4&\hspace{0.25cm}-X_9&\hspace{0.25cm}X_6-X_{10}&\hspace{0.25cm}0         &\hspace{0.25cm} X_9  \\
  X_{10}&\hspace{0.1cm} 0 &\hspace{0.25cm} 0  &\hspace{0.25cm}-X_3&\hspace{0.25cm} 0  &\hspace{0.25cm}0   &\hspace{0.25cm} X_5&\hspace{0.25cm}0   &\hspace{0.25cm}-X_7      &\hspace{0.25cm}-X_9      &\hspace{0.25cm} 0  \\
  \hline\hline
 \end{array}\end{eqnarray*}
\end{table}
\paragraph{Theorem 2.} {\em An optimal system of one-dimensional Lie subalgebras
NIB equation in the form  (\ref{eq:2}) is provided by those
generated by
\begin{eqnarray}
\begin{array}{r l r l }
1)  & A^1 = X_1,                 \hspace{2cm}    & 11) & A^{11} = \gamma_1\,X_6+X_{10}\\
2)  & A^2 = X_2,                 \hspace{2cm}    & 12) & A^{12} = \eta_1\,X_8+X_{10},\\
3)  & A^3 = X_3,                 \hspace{2cm}    & 13) & A^{13} = X_1+\eta_2\,X_8+X_{10},\\
4)  & A^4 = X_4,                 \hspace{2cm}    & 14) & A^{14} = X_2+\beta_1\,X_4+X_{10},\\
5)  & A^5 = X_5,                 \hspace{2cm}    & 15) & A^{15} = X_2+\gamma_2\,X_6+X_{10},\\
6)  & A^6 = X_6,                 \hspace{2cm}    & 16) & A^{16} = \alpha_1\,X_3+\gamma_3\,X_6+X_7,\\
7)  & A^7 = X_3+X_4,             \hspace{2cm}    & 17) & A^{17} = \alpha_2\,X_3+\gamma_4\,X_6+X_8,\\
8)  & A^8 = X_3+X_5,             \hspace{2cm}    & 18) & A^{18} = \gamma_5\,X_6+\eta_3\,X_8+X_9,\\
9)  & A^9 = X_3+X_6,              \hspace{2cm}   & 19) & A^{19} = X_4+\eta_4\,X_8+X_{10},\\
10) &A^{10}=X_4+X_{10}, \hspace{2cm}    & 20) & A^{20} =
\gamma_6\,X_6+\eta_5\,X_8+X_{10},
\end{array}\label{eq:30}
\end{eqnarray}
where $\alpha_i,\beta_1,\gamma_k,\eta_l \:(1\leq i\leq2,1\leq k\leq6,1\leq l\leq5)$ are arbitrary constants.}\\

\noindent{\it Proof:} Consider the symmetry algebra ${\cal
L}_{10}$ of Eq.~(\ref{eq:2}), whose adjoint representation was
determined in Table 4. Given a nonzero vector
\begin{eqnarray}
X=\sum_{i=1}^{10}\,a_i\,X_i,
\end{eqnarray}
we simplify $X$ by eliminating as many of the coefficients $a_i$
as possible by use of judicious applications of adjoint maps to
$X$. We perform the process through following cases:\\
\begin{table}
\centering{\caption{{\small Adjoint table for ${\cal
L}_{10}$}}}\label{table:4}\vspace{-0.5cm} {\tiny
\hspace{-2cm}\begin{eqnarray*}\begin{array}{l|c c c c c c c c c c}
\hline\hline
      \hspace{-0.2cm}&\hspace{-0.1cm} X_1 &\hspace{-0.4cm} X_2  &\hspace{-0.4cm}X_3       &\hspace{-0.3cm}X_4   &\hspace{-0.4cm}X_5          &\hspace{-0.2cm} X_6  &\hspace{-0.9cm}X_7        &\hspace{-0.6cm} X_8        &\hspace{-1cm}X_9       &\hspace{-0.3cm} X_{10}  \\ \hline
  X_1 \hspace{-0.2cm}&\hspace{-0.1cm} X_1 &\hspace{-0.4cm} X_2  &\hspace{-0.4cm}X_3       &\hspace{-0.3cm}X_4\!-\!sX_1 &\hspace{-0.4cm}X_5   &\hspace{-0.2cm} X_6  &\hspace{-0.9cm}X_7          &\hspace{-0.6cm} X_8        &\hspace{-1cm}X_9        &\hspace{-0.3cm} X_{10}  \\
  X_2 \hspace{-0.2cm}&\hspace{-0.1cm} X_1 &\hspace{-0.4cm} X_2  &\hspace{-0.4cm}X_3       &\hspace{-0.3cm} X_4  &\hspace{-0.4cm}X_5\!-\!sX_1 &\hspace{-0.2cm} X_6  &\hspace{-0.9cm}X_7\!-\!s\,X_2   &\hspace{-0.6cm} X_8        &\hspace{-1cm}X_9\!-\!sX_3       &\hspace{-0.3cm} X_{10}  \\
  X_3 \hspace{-0.2cm}&\hspace{-0.1cm} X_1 &\hspace{-0.4cm} X_2  &\hspace{-0.4cm}X_3       &\hspace{-0.3cm} X_4  &\hspace{-0.4cm}X_5        &\hspace{-0.2cm}X_6\!\!-\!\!sX_1&\hspace{-0.9cm}X_7      &\hspace{-0.6cm}X_8         &\hspace{-1cm}X_9         &\hspace{-0.3cm} X_{10}-sX_3  \\
  X_4 \hspace{-0.2cm}&\hspace{-0.1cm}e^sX_1&\hspace{-0.4cm} X_2   &\hspace{-0.4cm}X_3     &\hspace{-0.3cm}X_4   &\hspace{-0.4cm}X_5\!+\!sX_4  &\hspace{-0.2cm}X_6\!\!+\!\!sX_5\!\!-\!\!\frac{s^2}{2}X_4  &\hspace{-0.9cm}X_7         &\hspace{-0.6cm}X_8\!-\!sX_2  &\hspace{-1cm}X_9         &\hspace{-0.5cm} X_{10}  \\
  X_5 \hspace{-0.2cm}&\hspace{-0.1cm} X_1 &\hspace{-0.4cm}X_2\!\!+\!\!sX_1&\hspace{-0.4cm}X_3 &\hspace{-0.3cm}e^{-s}X_4 &\hspace{-0.4cm}X_5  &\hspace{-0.2cm} X_6  &\hspace{-0.9cm}X_7\!+\!(e^s\!\!-\!1)X_4&\hspace{-0.6cm}X_8\!+\!sX_5   &\hspace{-1cm}X_9         &\hspace{-0.3cm} X_{10}  \\
  X_6 \hspace{-0.2cm}&\hspace{-0.1cm} X_1 &\hspace{-0.4cm} X_2  &\hspace{-0.4cm}X_3\!+\!sX_1&\hspace{-0.3cm}X_4\!-\!sX_5 &\hspace{-0.4cm}X_5   &\hspace{-0.2cm} X_6  &\hspace{-0.9cm}X_7         &\hspace{-0.6cm} X_8         &\hspace{-1cm}X_9\!+\!sX_4\!-\!\frac{s^2}{2}X_5      &\hspace{-0.3cm}X_{10}+sX_5  \\
  X_7 \hspace{-0.2cm}&\hspace{-0.1cm} X_1 &\hspace{-0.4cm}e^sX_2  &\hspace{-0.4cm}X_3     &\hspace{-0.3cm} X_4   &\hspace{-0.4cm}X_5\!-\!sX_4 &\hspace{-0.2cm} X_6  &\hspace{-0.9cm}X_7         &\hspace{-0.6cm}X_8\!+\!sX_7   &\hspace{-1cm}e^{-s}\!X_9       &\hspace{-0.3cm} X_{10}  \\
  X_8 \hspace{-0.2cm}&\hspace{-0.1cm} X_1 &\hspace{-0.4cm} X_2  &\hspace{-0.4cm}X_3\!+\!sX_2&\hspace{-0.3cm} X_4 &\hspace{-0.4cm}e^{-s}X_5&\hspace{-0.2cm} X_6  &\hspace{-0.9cm}e^{-s}X_7 &\hspace{-0.6cm} X_8           &\hspace{-1cm}X_9\!\!-\!\!s(X_{10}\!\!-\!\!X_6)\!\!+\!\!(e^{-s}\!\!\!+\!\!s\!\!-\!\!1)X_7&\hspace{-0.3cm} X_{10}\!+\!(e^{-s}\!\!-\!1)X_7  \\
  X_9 \hspace{-0.2cm}&\hspace{-0.1cm} X_1 &\hspace{-0.4cm}X_2\!+\!sX_3&\hspace{-0.4cm}X_3 &\hspace{-0.3cm} X_4   &\hspace{-0.4cm}X_5        &\hspace{-0.2cm}X_6\!-\!sX_4&\hspace{-0.9cm}X_7\!+\!sX_9&\hspace{-0.6cm}X_8\!\!+\!\!s(X_6\!\!-\!\!X_{10})\!\!+\!\!\frac{s^2}{2}\!(X_9\!\!-\!\!X_4)&\hspace{-1cm}X_9         &\hspace{-0.3cm} X_{10}\!-\!sX_9  \\
X_{10}\hspace{-0.2cm}&\hspace{-0.1cm} X_1 &\hspace{-0.4cm} X_2 &\hspace{-0.4cm}e^s\,X_3   &\hspace{-0.3cm} X_4   &\hspace{-0.4cm}X_5        &\hspace{-0.2cm}X_6\!-\!sX_5&\hspace{-0.9cm}X_7     &\hspace{-0.6cm}X_8\!+\!sX_7   &\hspace{-1cm}e^s\!X_9        &\hspace{-0.5cm}X_{10}  \\
  \hline\hline
 \end{array}\end{eqnarray*}}
\end{table}

{\bf Case a.} Let $a_{10}\neq 0$. Scaling $X$ if necessary, we can
assume that $a_{10}=1$:
\begin{eqnarray}
X=\sum_{i=1}^{9}\,a_i\,X_i+X_{10}.
\end{eqnarray}
When $a_3=0$ or when $a_3\neq0$ by applying ${\rm
Ad}(\exp(a_3\,X_3)$ on $X$ to cancel the coefficient of $X_3$, one
can reduces $X$ to
\begin{eqnarray}
X'=a_1\,X_1+a_2\,X_2+a_4\,X_4+\cdots+X_{10}.
\end{eqnarray}
Now if $a_5=0$ or if $a_5\neq0$ by effecting ${\rm
Ad}(\exp(-a_5\,X_5)$ on $X'$, we can make the coefficient of $X_5$
vanish:
\begin{eqnarray}
X''=a_1\,X_1+a_2\,X_2+a_4\,X_4+a_6\,X_6+\cdots+X_{10}.
\end{eqnarray}
We can also cancel the coefficient of $X_7$. In the case which it
is nonzero, for $a_7>1$ by applying ${\rm Ad}(\exp(\ln
(a_7-1)\,X_7)$, for $a_7=1$ by applying ${\rm Ad}(\exp(-\ln
(2)\,X_7)$ and for $a_7<1$ by applying ${\rm Ad}(\exp(-\ln
(1-a_7)\,X_7)$ on $X''$ we change this coefficient to zero:
\begin{eqnarray}
X'''=a_1\,X_1+a_2\,X_2+a_4\,X_4+a_6\,X_6+a_8\,X_8+a_9\,X_9+X_{10}.
\end{eqnarray}
Furthermore, for different values of $a_9$, when it is either zero
or nonzero, this coefficient can be vanished; when $a_9\neq 0$ we
act ${\rm Ad}(\exp(a_9\,X_9)$ on $X'''$ to eliminate $a_9$.

{\bf Case a1.} Let $a_8\neq 0$ the additional action of ${\rm
Ad}(\exp(\frac{a_2}{a_8}\,X_4)$ change $X$ to the following form
which the coefficient of $a_2$ is deleted
\begin{eqnarray}
X''''=a_1\,X_1+a_4\,X_4+a_6\,X_6+a_8\,X_8+a_9\,X_9+X_{10}.
\end{eqnarray}

{\bf Case a1-1.} For $a_6\neq 0$ we effect ${\rm
Ad}(\exp(\frac{a_1}{a_6}\,X_3)$. Then the action of ${\rm
Ad}(\exp(\frac{a_4}{a_6}\,X_9)$ for $a_4\neq 0$ or when $a_4=0$ we
lead to the following form
\begin{eqnarray}
a_6\,X_6+a_8\,X_8+X_{10}.\label{eq:27}
\end{eqnarray}
Now, further simplification is impossible and every
one--dimensional subalgebra generated by an $X$ with
$a_{10},a_8,a_6\neq 0$ is equivalent to the subalgebra spanned by
(\ref{eq:27}). this vector field was introduced in part 20 of the
theorem.

{\bf Case a1-2.} If we change the condition of Case a1-1 to
$a_6=0$, the action of ${\rm Ad}(\exp(\frac{a_1}{a_4}\,X_1)$ for
$a_4\neq 0$ and then the action of ${\rm
Ad}(\exp(-\ln(\frac{1}{a_4})\,X_5)$ change $X''''$ to the below
final form (part 19 of the theorem)
\begin{eqnarray}
a_4\,X_4+a_8\,X_8+X_{10}.\label{eq:27}
\end{eqnarray}
In this case if $a_4=0$ we can effect ${\rm
Ad}(\exp(\ln(\frac{1}{a_1})\,X_4)$ for $a_1\neq0$ to reach to part
13 of the theorem, while $a_1=0$ suggests part 12.

{\bf Case a2.} Let $a_8$ in Case a1 be zero.

{\bf Case a2-1.} Moreover, if $a_6\neq0$ by effecting ${\rm
Ad}(\exp(\frac{a_1}{a_6}\,X_3)$ and then ${\rm
Ad}(\exp(\frac{a_4}{a_6}\,X_9)$ we tend to the form
\begin{eqnarray}
a_2\,X_2+a_6\,X_6+X_{10}.\label{eq:28}
\end{eqnarray}
In this form when $a_2\neq0$ the adjoint action of ${\rm
Ad}(\exp(\ln(\frac{1}{a_2})\,X_7)$ and when $a_2=0$ we find
sections 15 and 11 of the theorem resp.

{\bf Case a2-1.} Let $a_6=0$. Then the coefficient $a_1$ for
either $a_1\neq0$ by applying ${\rm
Ad}(\exp(\ln(\frac{1}{a_2})\,X_7)$ or $a_1=0$ will be zero:
\begin{eqnarray}
\tilde{X}=a_2\,X_2+a_4\,X_4+X_{10}.\label{eq:29}
\end{eqnarray}
Furthermore, the application of ${\rm
Ad}(\exp(\ln(\frac{1}{a_2})\,X_7)$ for $a_2\neq0$ introduces the
vector field $X_2+a_4\,X_4+x_{10}$, while for $a_2=0$ by effecting
${\rm Ad}(\exp(-\ln(\frac{1}{a_4})\,X_5)$ to change the
coefficient of $X_4$ equal to 1, we lead to the case $X_4+X_{10}$.

The remaining one--dimensional subalgebras are spanned by vectors
of  Case a, where $a_{10}=0$. Let $F^s _i(X)= {\rm
Ad}(\exp(s\,X_i)X)$ be a linear map, for $i = 1,\cdots,10$ and
every $X\in{\cal L}_{10}$. We continue the classification by a
similar method as above.

{\bf Case b.} Let $a_9\neq0$, we scale to make $a_9=1$. One can
set the coefficients of $X_3,X_2,X_5,X_7,X_1,X_4$ as zero (when
each of the coefficient is zero we eliminate it from the list) by
effecting
$F^{a_3}_2,F^{a_2}_4,F^{-a_5}_5,F^{-a_7}_7,F^{a_1}_3,F^{a_4}_9$
resp. Thus we tend to section 18.

{\bf Case c.} For $a_{10}=a_9=0, a_8\neq0$ and assuming $a_8=1$,
we can make the coefficients of $X_2,X_5,X_7, X_1,X_4$  vanish
(when those coefficients are not zero) by applying
$F^{a_2}_4,F^{-a_5}_5,F^{-a_5}_5,F^{-a_7}_7,F^{a_1}_1,F^{a_4}_9$
resp. This results in section 17.

{\bf Case d.} Let $a_{10}=a_9=a_8=0, a_7\neq0$ and scale if is
necessary to have $a_7=1$, we can cancel the coefficients of
$X_2,X_1,X_4,X_5$ by actions of
$F^{a_2}_2,F^{a_1}_3,F^{a_4}_9,F^{a_5}_{10}$ resp. to have section
16.

{\bf Case e.} Let $a_{10}=\cdots=a_7=0, a_6\neq0$. We can suppose
that $a_6=1$. Applying
$F^{a_1}_3,F^{a_4}_9,F^{a_5}_{10},F^{-a_2}_8$ on $X$ we can
changed the coefficients of $X_1,X_4,X_5,X_2$ to zero (if they are
not zero). Then if $a_3\neq0$ by acting ${\rm
Ad}(\exp(\ln(\frac{1}{a_3})\,X_{10})$ we find section 9 and when
$a_3=0$ we lead to section 6.

{\bf Case f.} Suppose that $a_{10}=\cdots=a_6=0, a_5\neq0$. By
assuming $a_5=1$, we can make the coefficients of $X_1,X_4,X_2$
vanish by $F^{a_1}_2,F^{-a_4}_4,F^{-a_2}_8$. Then if $a_3\neq0$ we
apply ${\rm Ad}(\exp(\ln(\frac{1}{a_3})\,X_6)$ to give section 8
and for $a_3=0$ we have section 5.

{\bf Case g.} Consider the case which $a_{10}=\cdots=a_5=0,
a_4\neq0$. We can assume that $a_4=1$. we can make the
coefficients of $X_1,X_2$ vanish by $F^{a_1}_1,F^{-a_2}_8$. In
addition for $a_3\neq0$, scaling the coefficient of $X_3$ by ${\rm
Ad}(\exp(\ln(\frac{1}{a_3})\,X_{10})$ we have section 7, while
$a_3=0$ leads to section 4.

{\bf Case h.} For $a_{10}=\cdots=a_4=0, a_3\neq0$ and assuming
$a_3=1$, the coefficients of $X_1,X_2$ will be vanished by the
actions of $F^{-a_1}_6,F^{-a_2}_8$ and $X$ reduces to section 3.

{\bf Case i.} If $a_{10}=\cdots=a_3=0, a_2\neq0$. Scaling $a_2=1$,
by effecting $F^{-a_1}_5$ the coefficient of $X_1$ will be
canceled which suggests section 2.

{\bf Case j.} Finally when $a_{10}=\cdots=a_2=0$ we have section 1
of the theorem.

There is not any more possible case for studying and the proof is
complete.\hfill\ $\diamondsuit$\\

Now, since $f,g$ in Eq.~(\ref{eq:2}) are functions of variables
$x,u$, so we project their derived optimal system to the space
$(x,u,f,g)$. The nonzero vector fields of $(x,u)-$projections of
(\ref{eq:30}) are as follows
\begin{eqnarray}\hspace{-2cm}\begin{array}{rl}
1)&Z^1=A^2=\partial_x,\\
2)&Z^2=A^3=u\,\partial_u , \\
3)&Z^3=A^4=f\,\partial_f + g\,\partial_g,\\
4)&Z^4=A^5=g[f\,\partial_f + g\,\partial_g], \\
5)&Z^5=A^6=f[f\,\partial_f + g\,\partial_g], \\
6)&Z^6=A^7=\partial_u -f\,\partial_f - g\,\partial_g, \\
7)&Z^7=A^8=\partial_u -g[f\,\partial_f + g\,\partial_g], \\
8)&Z^8=A^9=\partial_u -f[f\,\partial_f + g\,\partial_g],\\
9)&Z^9=A^{10}=u\,\partial_u -g\,\partial_g,\\
10)&Z^{10}=A^{11}=u\,\partial_u + (1-\gamma_1\,f)f\,\partial_f - f\,g\,\partial_g, \\
11)&Z^{11}=A^{12}=A^{13}=\eta_1\,u\,\partial_x + u\,\partial_u + f\,\partial_f +\eta_1\,g\,\partial_g\\
12)&Z^{12}=A^{14}=\partial_x + u\,\partial_u + (1-\beta_1)f\,\partial_f-\beta_1\,g\,\partial_g,\\
13)&Z^{13}=A^{15}=\partial_x + u\,\partial_u + f[(1-\gamma_2f)\partial_f-\gamma_2\,g\,\partial_g],\\
14)&Z^{14}=A^{16}=x\,\partial_x -\gamma_3\,f^2\,\partial_f+(1-\gamma_3\,f)\,g\,\partial_g,\\
15)&Z^{15}=A^{17}=u\,\partial_x + \alpha_2\,\partial_u -\gamma_4\,f^2\partial_f+(1-\gamma_4\,g)f\,\partial_g,\\
16)&Z^{16}=A^{18}=\eta_3\,u\,\partial_x + (x+\gamma_5\,u)\,\partial_u +(g-\gamma_5\,f^2)\partial_g,\\
17)&Z^{17}=A^{19}=\eta_4\,u\,\partial_x + u\,\partial_u +(\eta_4\,f-g)\partial_g,\\
18)&Z^{18}=A^{20}=\eta_5\,u\,\partial_x + u\,\partial_u
+f[(1-\gamma_6\,f)\,\partial_f+(\eta5-\gamma_6\,g)\,\partial_g].
\end{array}\label{eq:31}\end{eqnarray}
According to paper 7 of \cite{Ib3} we conclude that
\paragraph{Proposition 3.}{\em Let ${\cal L}_m:=\langle\,X_i: i = 1,
\cdots, m \,\rangle$ be an $m$--dimensional algebra. Denote by
$A^i\, (i = 1, \cdots, s,\, 0<s\leq m,\, s \in {\bf N})$ an
optimal system of one--dimensional subalgebras of ${\cal L}_m$ and
by $Z^i\, (i = 1, \cdots, t,\, 0<t\leq s,\, t\in {\bf N})$ the
projections of $A^i$, i.e., $Z^i = {\rm pr}(A^i)$. If equations
\begin{eqnarray}
g = g(x,u),\hspace{0.75cm} f =f(x,u),
\end{eqnarray}
are invariant with respect to the optimal system $Z^i$ then the
equation
\begin{eqnarray}
u_t +g(x,u)\,u_x=f(x,u),\label{eq:32}
\end{eqnarray}
admits the operators $Y^i=$ projection of $A^i$ on $(t, x, u)$.}

\paragraph{Proposition 4.} {\em Let Eq.~(\ref{eq:32}) and the equation
\begin{eqnarray}
u_t +\overline{g}(x,u)\,u_x=\overline{f}(x,u),\label{eq:33}
\end{eqnarray}
be constructed according to Proposition 3 via optimal systems
$Z^i$ and $\overline{Z}^i$ resp. If the subalgebras spanned on the
optimal systems $Z^i$ and $\overline{Z}^i$ resp. are similar in
${\cal L}_m$, then Eqs.~(\ref{eq:32}) and (\ref{eq:33}) are
equivalent with respect to the equivalence group
$G_m$ generated by ${\cal L}_m$.}\\

Regarding to Propositions 3 and 4, now our task is the
investigation for all non--equivalent equations in the form of
Eq.~(\ref{eq:2}) which admit ${\mathcal E}$--extensions of the
principal Lie algebra ${\mathcal L}_{{\mathcal E}}$, by one
dimension, that are equations of the form (\ref{eq:2}) such that
they admit, together with the one basic operator (\ref{eq:8-1}) of
${\mathcal L}_1$, also a second operator $X^{(2)}$. In each case
which this extension occurs, we indicate the corresponding
coefficients $f , g$ and the additional operator $X^{(2)}$.

By using the algorithm for operators $Z^i\,(i=1,\cdots,18)$ to
$f,g$ one can find vector fields $X^{(2)}$~s. Let we consider the following examples.\\
\mbox{ }\hspace{3mm}   We select the operator
\begin{eqnarray}
Z^5=f^2\,\frac{\partial}{\partial
f}+f\,g\,\frac{\partial}{\partial g}.
\end{eqnarray}
The characteristic equations are
\begin{eqnarray}
\frac{df}{f^2}=\frac{dg}{f\,g},
\end{eqnarray}
so the invariants of $Z^5$ are
\begin{eqnarray}
I_1=x, \hspace{0.5cm} I_2= u, \hspace{0.5cm}
I_3=\frac{f}{g}.\label{eq:35}
\end{eqnarray}
In this case there are no invariant equations because the
necessary condition for existence of invariant solutions is not
fulfilled (see \cite{Ov}, Section 19.3), that is, invariants
(\ref{eq:35}) cannot be solved with respect to $f$ and $g$ since
$I_2$ is not an invariant function of $I_1$. Similar results are
hold for vector fields $Z^3$ and $Z^4$.

\begin{table} \centering{\caption{The result of the
classification}}\label{table:5} \vspace{-0.35cm}
{\small\begin{eqnarray*} \hspace{-0.75cm}\begin{array}{l l l l l
l} \hline\hline
  N  &\hspace{0.0cm} Z      &\hspace{0.1cm} \mbox{Invariant $\lambda$}  &\hspace{0.1cm}  \mbox{Equation} &\hspace{-1cm} \mbox{ Additional operator}\,X^{(2)} \\
  \hline
  & &&& \\[-3mm]
  1  &\hspace{0.0cm} Z^1    &\hspace{0.1cm} u   &\hspace{0.1cm}  u_t\!+\!\Psi u_x\!=\!\Phi               &\hspace{0.1cm} \frac{\partial}{\partial x} \\[2mm]
  2  &\hspace{0.0cm} Z^2    &\hspace{0.1cm} x   &\hspace{0.1cm}  u_t\!+\!\Psi u_x\!=\!\Phi               &\hspace{0.1cm} \frac{\partial}{\partial u} \\[2mm]
  3  &\hspace{0.0cm} Z^6    &\hspace{0.1cm} x   &\hspace{0.1cm}  u_t\!+\!e^{\Psi\!-\!u} u_x\!=\!e^{\Phi\!-\!u}     &\hspace{0.1cm} t\,\frac{\partial}{\partial t}+\frac{\partial}{\partial u}\\[2mm]
  4  &\hspace{0.0cm} Z^7    &\hspace{0.1cm} x   &\hspace{0.1cm}  u_t\!+\![\Psi\!+\!1/u]u_x\!=\![\Psi\!+\!1/u]\Phi  &\hspace{0.1cm} x\,\frac{\partial}{\partial t}+\frac{\partial}{\partial u}\\[2mm]
  5  &\hspace{0cm} Z^8    &\hspace{0.1cm}x    &\hspace{0.1cm}  u_t\!+\![\Phi\!+\!1/u]\Psi u_x\!=\!\Phi\!+\!1/u   &\hspace{0.1cm} u\,\frac{\partial}{\partial t}+\frac{\partial}{\partial u}\\[2mm]
  6  &\hspace{0cm} Z^9    &\hspace{0.1cm} x   &\hspace{0.1cm}  u_t\!+\!u\Psi u_x\!=\!\Phi   &\hspace{0.1cm}  t\,\frac{\partial}{\partial t}+u\,\frac{\partial}{\partial u}\\[2mm]
  7  &\hspace{0cm} Z^{10} &\hspace{0.1cm} x   &\hspace{0.1cm}  u_t\!\pm\!(\gamma_1 f\!-\!1)\Psi u_x\!=\!\frac{u\Phi}{1\!+\!\gamma_1 u \Phi} &\hspace{0.1cm} u[\gamma_1\frac{\partial}{\partial t}+\frac{\partial}{\partial u}]\\[2mm]
  8  &\hspace{0cm}Z^{11}  &\hspace{0.1cm}x-\eta_1 u &\hspace{0.1cm}u_t\!+\![\Psi\!-\!\eta_1u\Phi]u_x\!=\!u\Phi   &\hspace{0.1cm} u[\eta_1\frac{\partial}{\partial x}\!+\!\frac{\partial}{\partial u}\!],\\[2mm]
  & &&& \hspace{0.1cm}  \!\frac{\partial}{\partial t}\!+\!u[\eta_1\frac{\partial}{\partial x}\!+\!\frac{\partial}{\partial
  u}\!]\\[2mm]
   9  &\hspace{0cm}Z^{12} &\hspace{0.1cm}\ln u-x &\hspace{0.1cm} u_t\!+\!u^{\beta_1}\!\Psi u_x\!=\!u^{1\!-\!\beta_1}\Phi &\hspace{0.1cm} \beta_1 t\frac{\partial}{\partial t}+\frac{\partial}{\partial x}+u\,\frac{\partial}{\partial u} \\[2mm]
  10 &\hspace{0cm}Z^{13}(\gamma_2\neq0)&\hspace{0.1cm}\ln u\!-\!x  &\hspace{0.1cm}u_t\!+\![(\frac{\gamma_2u\Phi}{\gamma_2u\Phi-1}-1)\Psi]^{\frac{1}{\gamma_2}}u_x\!=\!\frac{\gamma_2u\Phi}{\gamma_2u\Phi\!-\!1} &\hspace{0.1cm}\gamma_2u\frac{\partial}{\partial t}+\frac{\partial}{\partial x}+u\frac{\partial}{\partial u}\\[2mm]
  11 &\hspace{0cm}Z^{13}(\gamma_2=0)&\hspace{0.1cm}\ln u\!-\!x   &\hspace{0.1cm}u_t\!+\!\Psi u_x\!=\!u\Phi  &\hspace{0.1cm} \frac{\partial}{\partial x}\\[2mm]
  12 &\hspace{0cm}Z^{14}(\gamma_3\neq0) &\hspace{0.1cm} u   &\hspace{0.1cm}u_t\!+\!x^{1\!-\!\gamma_3\![\Phi\!+\!\frac{1}{\gamma_3\ln\! x}]}\!\Psi u_x\!=\!\Phi\!+\!\frac{1}{\gamma_3\ln x}   &\hspace{0.1cm}\gamma_3\,u\,\,\frac{\partial}{\partial t}+x\,\frac{\partial}{\partial x}\\[2mm]
  13 &\hspace{0cm}Z^{14}(\gamma=0) &\hspace{0.1cm} u  &\hspace{0.1cm}u_t\!+\!x\Psi u_x\!=\!\Phi   &\hspace{0.1cm} x\,\frac{\partial}{\partial x} \\[2mm]
  14 &\hspace{0cm}Z^{15}(\alpha_2=\gamma_4=0) &\hspace{0.1cm}x   &\hspace{0.1cm}u_t\!+\![\Psi\!-\!\frac{x}{u}]\Phi u_x\!=\!\Phi   &\hspace{0.1cm} u\,\frac{\partial}{\partial x}\\[2mm]
  15 &\hspace{0cm}Z^{15}(\alpha_2\neq0,\gamma_4=0) &\hspace{0.1cm}\frac{u}{\alpha_2}-x   &\hspace{0.1cm} u_t\!+\![\Psi\!-\!\frac{x}{u}]\Phi u_x\!=\!\Phi    &\hspace{0.1cm} u\,\frac{\partial}{\partial x}+\alpha_2\frac{\partial}{\partial u}\\[2mm]
  16 &\hspace{0cm}Z^{15}(\alpha_2=0,\gamma_4\neq0) &\hspace{0.1cm}u   &\hspace{0.1cm}u_t\!+\![\frac{1}{\gamma_4}\!\pm\!(\frac{1}{\gamma_4}\Phi\!-\!\frac{u}{x})\Psi] u_x\!=\!\Phi\!-\!\frac{\gamma_4 u}{x}   &\hspace{0.1cm} \gamma_4 u \frac{\partial}{\partial t}+u\,\frac{\partial}{\partial x}\\[2mm]
  17 &\hspace{0cm}Z^{15}(\alpha_2\neq0,\gamma_4\neq0)   &\hspace{0.1cm}\frac{u}{\alpha_2}-x   &\hspace{0.1cm}u_t\!+\!\frac{1}{\gamma_4}[1\!\pm\!\Phi\Psi]u_x\!=\!\Phi  &\hspace{0.1cm} \gamma_4 u \frac{\partial}{\partial t}+u\,\frac{\partial}{\partial x}+\alpha_2\frac{\partial}{\partial u}\\[2mm]
  18 &\hspace{0cm}Z^{16}(\eta_3=\gamma_5=0) &\hspace{0.1cm}x   &\hspace{0.1cm}u_t\!+\![\Psi\!+\!e^{u/x}]u_x\!=\!\Phi   &\hspace{0.1cm} x\frac{\partial}{\partial u} \\[2mm]
  19 &\hspace{0cm}Z^{16}(\eta_3\neq0,\gamma_5=0)  &\hspace{0.1cm}\frac{1}{2}(\eta_3u^2\!-\!x^2)   &\hspace{0.1cm}u_t\!+\![\Psi\!+\!e^{u/x}]u_x\!=\!\Phi  &\hspace{0.1cm} \eta_3 u \frac{\partial}{\partial x}+x\frac{\partial}{\partial u} \\[2mm]
  20 &\hspace{0cm}Z^{16}(\eta_3=0,\gamma_5\neq0)  &\hspace{0.1cm}x   &\hspace{0.1cm}u_t\!+\![\gamma_5\Phi^2\!\pm\!(\gamma_5u\!+\!x)^{\frac{1}{\gamma_5}}\!\Psi]u_x\!=\!\Phi     &\hspace{0.1cm} [\gamma_5 u +x]\frac{\partial}{\partial u} \\[2mm]
  21 &\hspace{0cm}Z^{16}(\eta_3\neq0,\gamma_5\neq0)&\hspace{0.1cm}B  &\hspace{0.1cm}u_t\!+\![\gamma_5\Phi^2\!\pm\! e^{x/(\eta_3 u)}\!+\!\Psi]u_x\!=\!\Phi    &\hspace{0.1cm} \eta_3 u \frac{\partial}{\partial x}+[\gamma_5 u +x]\frac{\partial}{\partial u} \\[2mm]
  22 &\hspace{0cm}Z^{17}(\eta_4\neq0) &\hspace{0.1cm}\frac{x}{\eta_4}\!-\!u &\hspace{0.1cm}u_t\!+\![\eta_4\Phi\!\pm\! u\Psi]u_x\!=\!\Phi    &\hspace{0.1cm} t \frac{\partial}{\partial t} + u[\eta_4\frac{\partial}{\partial x}+\frac{\partial}{\partial u}]\\[2mm]
  23 &\hspace{0cm}Z^{17}(\eta_4=0)    &\hspace{0.1cm}x   &\hspace{0.1cm}u_t\!+\!u\Psi u_x\!=\!\Phi     &\hspace{0.1cm} t\frac{\partial}{\partial t} + u\frac{\partial}{\partial u}\\[2mm]
  24 &\hspace{0cm}Z^{18}(\eta_5=\gamma_6=0) &\hspace{0.1cm}u\!-\!\frac{x}{\eta_5}   &\hspace{0.1cm}u_t\!+\!\Psi u_x\!=\!u\Phi     &\hspace{0.1cm}  u\frac{\partial}{\partial u}\\[2mm]
  25 &\hspace{0cm}Z^{18}(\eta_5\neq0,\gamma_6=0)      &\hspace{0.1cm} \ln u\pm x     &\hspace{0.1cm}u_t\!+\!\eta_5[u\Phi\!+\!\Psi]u_x\!=\!u\Phi    &\hspace{0.1cm}   \eta_5u\frac{\partial}{\partial x}+u\frac{\partial}{\partial u}\\[2mm]
  26 &\hspace{0cm}Z^{18}(\eta_5=0,\gamma_6\neq0)      &\hspace{0.1cm}x     &\hspace{0.1cm}u_t\!\pm\![\frac{u\Phi}{1\!+\!\gamma_6u\Phi}\!-\!\frac{1}{\gamma_6}]\Psi u_x\!=\!\frac{u\Phi}{1\!+\!\gamma_6u \Phi}     &\hspace{0.1cm}  \gamma_6u\frac{\partial}{\partial t}+u\frac{\partial}{\partial u} \\[2mm]
  27 &\hspace{0cm}Z^{18}(\eta_5\neq0,\gamma_6\neq0)   &\hspace{0.1cm}u\!-\!\frac{x}{\eta_5} &\hspace{0.1cm}u_t\!+\![\frac{\eta_5}{\gamma_6}\!\pm\!(\frac{u\Phi}{1\!+\!\gamma_6u\Phi}\!-\!\frac{1}{\gamma_6})\Psi] u_x\!=\!\frac{u\Phi}{1\!+\!\gamma_6u\Phi}     &\hspace{0.1cm} \gamma_6u\frac{\partial}{\partial t}+\eta_5u\frac{\partial}{\partial x}+u\frac{\partial}{\partial u} \\[2mm]
  \hline\hline
\end{array}\end{eqnarray*}}
\end{table}

Let perform the algorithm for another example, by considering
vector filed
\begin{eqnarray}
Z^{14}=x\,\frac{\partial}{\partial x}
-\gamma_3\,f^2\,\frac{\partial}{\partial
f}+(1-\gamma_3f)g\,\frac{\partial}{\partial g},
\end{eqnarray}
then the characteristic equations corresponding to $Z^{14}$ is
\begin{eqnarray}
\frac{dx}{x}=\frac{df}{-\gamma_3\,f^2}=\frac{dg}{(1-\gamma_3f)g},
\end{eqnarray}
which determines invariants for $\gamma_3\neq0$. Invariants can be
taken in the following form
\begin{eqnarray}
I_1=u, \hspace{0.5cm} I_2= f-\frac{1}{\gamma_3\,\ln x},
\hspace{0.5cm} I_3=\frac{g}{x^{1-\gamma_3\,f}}. \label{eq:34}
\end{eqnarray}
From the invariance equations we can write
\begin{eqnarray}
I_2=\Phi(I_1), \hspace{1cm} I_3= \Psi(I_1).
\end{eqnarray}
It results in the forms
\begin{eqnarray}
f=\Phi(\lambda)+\frac{1}{\gamma_3\,\ln x},\hspace{1cm}
g=x^{1-\gamma_3[\Phi(\lambda)+\frac{1}{\gamma_3\,\ln x}]}
\,\Psi(\lambda),
\end{eqnarray}
where $\lambda=u$. For the case which $\gamma_3=0$, $Z_{14}$ is in
the form $x\,\partial_x +g\,\partial_g$ and similar way shows that
$f=\Phi(\lambda)$ and $g=x\,\Psi(\lambda)$ when $\lambda=u$.

From Proposition 3 applied to the operator $Z^{14}$ we obtain the
following additional operators for the introduced functions $f$
and $g$ as above
\begin{eqnarray}
X^{(2)}= \gamma_3\,u\,\frac{\partial}{\partial t} +
x\,\frac{\partial}{\partial x},
\end{eqnarray}
When $\gamma_3\neq0$, while for $\gamma_3=0$ we find $X^{(2)}=
x\,\partial_x$ as an additional operator.

Repeating the algorithm for other $Z^i$~s of (\ref{eq:31}), the
preliminary group classification of an NIB equation of the form
(\ref{eq:2}) admitting an extension ${\mathcal L}_2$ of the
principal Lie algebra ${\mathcal L}_1$ is listed in Table 5. In
this table we assume that
\begin{eqnarray}
B = x + \gamma_3u + \frac{\eta_1\,u}{\gamma_3}\,\,{\rm
LambertW}(A)
\end{eqnarray}
where
\begin{eqnarray}
A=-\frac{\gamma_3}{\eta_1u}\,e^{-\gamma_3^2/\eta_1},
\end{eqnarray}
and $\Phi=\Phi(\lambda), \Psi=\Psi(\lambda)$ are arbitrary
functions of the invariant $\lambda$.
%
%%%%%%%%%%%%%%%%%%%%%%%%%%%%%%%%%%%%%%%%%%%%%%%%%%%%%%%%%%%%%%
\section{Conclusion}
A symmetry analysis for generalized non--homogeneous inviscid
Burger's equations of class (\ref{eq:2}) led to find the structure
of point infinitesimal generators as well as equivalence
operators. In addition we found the modified forms of point and
projective symmetries of IBE in the paper \cite{Na} among with the
general form of point infinitesimals of a generalized case with
respect to IBE and results of \cite{SM}. Also, equivalence
classification is given of the equation admitting an extension by
one of the principal Lie algebra of the equation. The paper is one
of few applications of a new algebraic approach to the problem of
group classification: the method of preliminary group
classification. Derived results are summarized in Table 5.

%%%%%%%%%%%%%%%%%%%%%%%%%%%%%%%%%%%%%%%%%%%%%%%%%%%%%%%%%%%%%%%%%%%%%%%%%%%%%%%%%%%%%%%%%%%%%%%%%%%%%%%%%%%%

\end{document}